\documentclass[11pt]{article}
\usepackage{mathrsfs, cite}
\usepackage{amsmath,amsfonts,amssymb,rotating,amsthm}
\usepackage{psfrag,eepic,color}
\usepackage{framed}

\def\qed{\nopagebreak\hfill{\rule{4pt}{7pt}}}
\def\proof{\noindent {\it{Proof.} \hskip 2pt}}
\newtheorem{theo}{Theorem}[section]

\newtheorem{lem}[theo]{Lemma}

\newtheorem{ex}{Example}[section]
\theoremstyle{remark}

\begin{document}
\begin{center}
{\Large \bf $\ell$-Log-momotonic and Laguerre Inequality of P-recursive Sequences}
\end{center}

\begin{center}
{\large Guo-Jie Li}\\[9pt]
School of Science\\
Hainan University\\
Haikou 570208, China\\[6pt]
{\tt gj\_li@hainanu.edu.cn}
\end{center}

\vspace{0.3cm} \noindent{\bf Abstract.}

We consider $\ell$-log-momotonic sequences  and Laguerre inequality of order two for sequences $\{a_n\}_{n \ge 0}$ such that
\[
\frac{a_{n-1}a_{n+1}}{a_n^2}  = 1 + \sum_{i=1}^m  \frac{r_i(\log n)}{n^{\alpha_i}} + o\left( \frac{1}{n^{\beta}} \right),
\]
where $m$ is a nonnegative integer, $\alpha_i$ are real numbers, $r_i(x)$ are rational functions of $x$ and
\[
0 < \alpha_1 < \alpha_2 < \cdots < \alpha_m < \beta.
\]
We will give a sufficient condition on $\ell$-log-momotonic sequences  and Laguerre inequality of order two for $n$ sufficiently large.
 Many P-recursive sequences fall in this frame. At last,  we will give a method to find the  $N$ such that for any $n\geq N$,  log-momotonic inequality of order three and Laguerre inequality of order two holds.

\vskip 15pt


\vskip 15pt


\section{Introduction}
The Tur\'{a}n inequalities  arise in the study
of Maclaurin coefficients of real entire functions in the Laguerre-P\'{o}lya class \cite{Sze}.
A sequence $\{a_n\}_{n\geq 0}$ of real numbers is said to satisfy the Tur\'{a}n inequalities
or to be log-concave, if
\[
a_n^2 -a_{n-1} a_{n+1} \geq 0, \quad \forall\, n\geq 1.
\]
The Tur\'an inequalities are also called Newton's inequality \cite{Wag77, Nic}. For more results on the log-concavity, we refer to \cite{Des, McN, Sta89}.

Let $\varphi$ be the operator defined on sequences $\{a_n\}_{n\geq0}$ by
\[
\varphi\{a_n\}_{n\geq0}=\{b_n\}_{n\geq0},
\]
where $b_n=a_{n+1}/a_n$. A sequence $\{a_n\}_{n\geq0}$ is log-monotonic of order $k$ if for $r$ odd and not greater than $k-1$, the sequence
$\varphi^r\{a_n\}_{n\geq0}$ is log-concave and for $r$ even and not great than $k-1$, the sequence $\varphi^r\{a_n\}_{n\geq0}$ is log-convex.
A sequence $\{a_n\}_{n\geq0}$ is called infinitely log-monotonic if it is log-monotonic of order $k$ for all integers $k\geq1$.
Log-monotonic of order $\ell$ is  also called as $\ell$-log-monotonic.

A sequence $\{a_n\}_{n\geq0}$ is said to be ratio log-concave if $\{a_{n+1}/a_n\}_{n\geq0}$ is log-concave. Chen, Guo and Wang\cite{CGW}
proved that under certain initial condition, the ratio log-concavity of a sequence $\{a_n\}_{n\geq k}$ of positive numbers implies that the sequence
$\{\sqrt[n]{a_n}\}_{n\geq k}$ is strictly log-concave.

A function $f$ is said to be completely monotonic on an interval $I$ if $f$ has derivatives of all orders on $I$ and
\[
(-1)^n f^{(n)}(x) \geq 0
\]
for $x \in I$ and all integers $n\geq0$. Chen, Guo and Wang\cite{CGW} established a connection between completely momotonic functions and infinitely
 log-monotonic sequences and showed the sequences of the Bernoulli numbers, the Catalan numbers and the central binomial coefficients are infinitely log-monotonic.

The Laguerre inequality\cite{Larson} arises in the study of the real polynomials with only real zeros and the Laguerre-P\'olya class consisting of real entire functions. Recall
that a  real entire function
\[
\varphi(x)  =\sum_{k=0}^{\infty}  \gamma_k  \frac{x^k}{k!}
\]
is said to be in the Laguerre-P\'olya class, denoted by $\varphi(x) \in \mathcal{LP}$, if it can be represented in the form
\[
\varphi(x) =c x^m e^{-\alpha x^2+\beta x} \prod_{k=1}^{\infty} (1+x/x_k) e^{-x/x_k},
\]
where $c$, $\beta$, $x_k$ are real numbers, $\alpha\geq0$, $m$ is a nonnegative integer and $\sum x_k^{-2}<\infty$. For more background on the theory of the
$\mathcal{LP}$ class, we refer to \cite{Levin} and \cite{Rahman}.

Recall that if a polynomial $f(x)$ satisfies
\[
f'(x)^2-f(x)f''(x)\geq 0,
\]
then it is called to satisfy Laguerre inequality. Laguerre stated that if $f(x)$ is a polynomial with only real zeros, then the Laguerre inequality holds for $f(x)$.
Gasper\cite{Gasper} used it as an important tool to deal with the positivity of special function.

In 1913, Jensen\cite{Jensen} found a $n$-th generalization of the Laguerre inequality
\[
L_n(f(x)): =\frac{1}{2} \sum_{k=0}^{2n} (-1)^{n+k} \binom{2n}{k}  f^k(x)  f^{2n-k}(x) \geq 0,
\]
where $f^k(x)$ denotes the $k$-th derivative of $f(x)$. It yields the classical Laguerre inequality for $n=1$.

Yang and Wang\cite{YangW} defined a sequence $\{a_n\}_{n\geq0}$ satisfies Laguerre inequality of order $m$ if
\[
L_m(a_n): =\frac{1}{2} \sum_{k=0}^{2m}  (-1)^{k+m}  \binom{2m}{k}  a_{n+k}  a_{2m-k+n} \geq 0.
\]
For $m=1$, the above inequality reduces to log-concavity of the sequence $\{a_n\}$. Yang and Wang\cite{YangW} proved some combinatorial sequences satisfy the Laguerre inequality of order two. Wagner\cite{Wagner} recently proved that the partition function $p(n)$ satisfied all of the Laguerre inequalities of order $m$ as $n\to \infty$ and proposed a conjecture
on the thresholds of the $m$-th Laguerre inequalities of $p(n)$ for $m\leq 10$.

We focus on the behaviour of a sequence $\{a_n \}_{n \ge 0}$ when $n$ is sufficiently large.  We say a sequence satisfies Log-momotonic inequality and Laguerre inequality asymptotically if for $n$ sufficiently large the  Log-monotonic inequality and Laguerre inequality holds.  We aim to give a criterion on the $\ell$-log-momotonic inequality  and Laguerre inequality of order two for P-recursive sequences.   A P-recursive sequence of order $d$ satisfies a recurrence relation of the form
\[
a_n=r_1(n)a_{n-1}+r_2(n)a_{n-1}+\cdots+r_d(n)a_{n-d},
\]
where $r_i(n)$ are rational functions of $n$, see \cite[Section 6.4]{Sta}.

By the asymptotic estimation given by Birkhoff and Trjitzinsky \cite{Birk} and developed by Wimp and Zeilberger \cite{Wim}, Hou and Li\cite[Theorem 4.1]{Houli} proved that many P-recursive sequences $\{a_n\}_{n \ge 0}$ satisfy
\begin{equation}\label{un}
u_n = \frac{a_{n-1}a_{n+1}}{a_n^2} = 1 + \sum_{i=1}^m  \frac{r_i(\log n)}{n^{\alpha_i}} + o\left( \frac{1}{n^{\beta}} \right),
\end{equation}
where $m$ is a nonnegative integer, $\alpha_i$ are real numbers, $r_i(x)$ are rational functions of $x$ and
\[
0 < \alpha_1 < \alpha_2 < \cdots < \alpha_m < \beta.
\]
Hou and Li\cite[Lemma 2.1]{Houli} also give a estimation we will use later.
\begin{lem}\label{dif-log}
Let $r(x)$ be a rational function of $x$ and $K$ be a positive integer. Then
\[
r(\log (n+1)) - r( \log n)
=  \sum_{i=1}^K \frac{r_i(\log n)}{n^i} + o\left( \frac{1}{n^K} \right),
\]
for some rational functions $r_i(x)$. Similarly,
\[
r(\log (n-1)) - r( \log n)
=  \sum_{i=1}^K \frac{\tilde{r}_i(\log n)}{n^i} + o\left( \frac{1}{n^K} \right),
\]
for some rational functions $\tilde{r}_i(x)$. If $r(x)$ is a polynomial, are the $r_i(x), \tilde{r}_i(x)$'s. Moreover,
\[
r_1(x) = - \tilde{r}_1(x) = r'(x) \quad \mbox{and} \quad
r_2(x)=\tilde{r}_2(x)= \frac{r''(x)-r'(x)}{2}.
\]
\end{lem}
With the asymptotic form \eqref{un} and Lemma~\ref{dif-log}, we are able to give a sufficient condition on the $\ell$-log-monotonic inequality  in Section~\ref{sec2}. Then in Section~\ref{sec3}, we will give a sufficient condition on the asymptotic Laguerre inequality of order two.
We then  give a method to find the $N$ such that for any $n\geq N$  the  log-monotonic inequality of order three holds in section~\ref{sec4}. At last, we give a method
to find the $N$ such that for any $n\geq N$ the Laguerre inequality of order two holds in section ~\ref{sec5}.

\section{Asymptotically $\ell$-log-monotonic }\label{sec2}
In this section, we  give a sufficient condition such that the sequence is $\ell$-log-momotonic.

We first give a lemma we will use later.

\begin{lem}\label{lem-2.1}

Let $r(x)$ be a rational function of $x$, $\xi(n)=\frac{r(\log n)}{n^\gamma}$, $\gamma$, $\alpha$ are real numbers and $\gamma>0$, $\alpha>0$, then we have the following estimation
\[
(n+1)^\alpha \xi(n)-n^\alpha \xi(n+1)=r(\log n) n^{\alpha-\gamma-1} \left(\alpha+\gamma+o(1)\right),
\]
and
\begin{align*}
&\xi(n+1) (n-1)^\alpha  n^\alpha + \xi(n-1) (n+1)^\alpha  n^\alpha - 2\xi(n) (n-1)^\alpha (n+1)^\alpha\\
=&r(\log n) n^{2\alpha-\gamma-2} \big( \left (\alpha+\gamma \right) \left( \alpha+\gamma+1 \right)+o(1) \big),
\end{align*}
and
\begin{align*}
&\xi(n-1) \xi(n+1) n^{2\alpha} - \xi^2(n) (n+1)^\alpha (n-1)^\alpha\\
=&r^2(\log n) n^{2\alpha-2\gamma-2} \left( \alpha+\gamma+o(1) \right).
\end{align*}

\end{lem}
\proof Since $n \to \infty$,
\[
\lim_{n\to \infty}\frac{r'(\log n)}{r(\log n)}=0.
\]
By Lemma ~\ref{dif-log} we have
\[
\frac{r(\log(n+1))}{r(\log n)} = 1 + \frac{r'(\log n)}{ r(\log n) n} + o \left( \frac{1}{n^2} \right)
\]
and
\[
\frac{r(\log(n-1))}{r(\log n)} = 1-  \frac{r'(\log n)}{ r(\log n) n} + o \left( \frac{1}{n^2} \right).
\]
Then we can deduce that
\begin{align*}
&(n+1)^\alpha \xi(n) - n^\alpha \xi(n+1)\\
=&\frac{r(\log n) n^\alpha}{n^\gamma} \big[ \left( 1+\frac{1}{n} \right)^\alpha - \frac{r(\log(n+1))}{r(\log n)} \left( 1+\frac{1}{n} \right)^{-\gamma} \big]\\
=&r(\log n) n^{\alpha-\gamma-1} \left[ \alpha+\gamma+o(1) \right].
\end{align*}
And
\begin{align*}
&\xi(n+1) (n-1)^\alpha  n^\alpha + \xi(n-1) (n+1)^\alpha  n^\alpha - 2\xi(n) (n-1)^\alpha (n+1)^\alpha\\
=&r(\log n) n^{2\alpha-\gamma} \big[ \frac{r(\log(n+1))}{r(\log n)} \left( 1-\frac{1}{n} \right)^\alpha \left( 1+\frac{1}{n} \right)^{-\gamma} + \\
 & \qquad \qquad \qquad  \frac{r(\log(n-1))}{r(\log n)} \left( 1+\frac{1}{n} \right)^\alpha  \left( 1-\frac{1}{n} \right)^{-\gamma}  -  2\left( 1-\frac{1}{n^2} \right)^\alpha \big] \\
=&r(\log n) n^{2\alpha-\gamma-2} \left[ \alpha(\alpha+1) + \gamma(\gamma+1) + 2\alpha\gamma + o(1) \right]\\
=&r(\log n) n^{2\alpha-\gamma-2} \left[ (\alpha+\gamma) (\alpha+\gamma+1) + o(1) \right].
\end{align*}
And
\begin{align*}
&\xi(n-1) \xi(n+1) n^{2\alpha} - \xi^2(n) (n+1)^\alpha (n-1)^\alpha\\
=&r^2(\log n) n^{2\alpha-2\gamma} \big[ \frac{r(\log(n+1)r(\log(n-1)))}{r^2(\log n)}  \left( 1-\frac{1}{n^2} \right)^{-\gamma} - \left( 1-\frac{1}{n^2} \right)^\alpha \big]\\
=&r^2(\log n) n^{2\alpha-2\gamma-2} \left[ \alpha+\gamma+o(1) \right].
\end{align*}

 Now I will give the main Theorem of this section.
\begin{theo}\label{Th-2.3}
Let $\{a_n\}_{n \ge 0}$ be a sequence such that \eqref{un} holds.
Assume that $\alpha_m- \alpha_1 \ge 1$ and denote $\ell=\lfloor \alpha_m/\alpha_1\rfloor$. If $r_1(x)>0$ for $x$ sufficiently large,
then $\{a_n\}_{n \ge 0}$ satisfies the  $\ell$-log-monotonic inequality  asymptotically.
\end{theo}
\proof Let
\[
u_n = \frac{a_{n-1}a_{n+1}}{a_n^2}, \quad
b_n = \frac{a_{n+1}}{a_n}.
\]
Denote
\[
\alpha_1 = \alpha,  \quad
u_n = 1 + \frac{\xi(n)}{n^{\alpha}}.
\]
Then $\xi(n)$ has the form
\[
\xi(n)= \sum_{i=1}^m  \frac{r_i(\log n)}{n^{\alpha_i-\alpha}} + o\left( \frac{1}{n^{\beta-\alpha}} \right).
\]
Let $\varphi$ be the operator defined on  sequence $\{a_n\}$ by
\[
\varphi\{a_n\} = \{b_n\} = \{\frac{a_{n+1}}{a_n}\},
\]
It is easy to check that
\[
\frac{\varphi\{a_{n-1}\} \varphi\{a_{n+1}\} }{ (\varphi\{a_{n}\})^2 } = \frac{ b_{n-1} b_{n+1} }{b_n^2} = \frac{u_{n+1}}{u_n}.
\]
By Lemma ~\ref{lem-2.1}, we have
\begin{align*}
&u_n-u_{n+1}\\
=&\frac{(n+1)^\alpha \xi(n) - n^\alpha \xi(n+1)}{n^\alpha (n+1)^\alpha}\\
=&\frac{r_1(\log n)}{n^{\alpha+1}} \left(\alpha+o(1) \right) > 0, \quad n \to \infty.
\end{align*}
Then we have
\[
\frac{\varphi\{a_{n-1}\} \varphi\{a_{n+1}\} }{ (\varphi\{a_{n}\})^2 } = 1 + \frac{u_{n+1}-u_n}{u_n} = 1+(-1) \frac{r_1(\log n)}{n^{\alpha+1}} \left( \alpha+o(1) \right).
\]
Repeating this argument, we have
\[
\frac{\varphi^k\{a_{n-1}\} \varphi^k\{a_{n+1}\} }{ (\varphi^k\{a_{n}\})^2 } = 1 + (-1)^k \frac{r_1(\log n)}{n^{\alpha+k}} \left( (\alpha)_k+o(1) \right),
\]
where
\begin{quote}
$(\alpha)_k=1$, \quad $k=0$. \\
$(\alpha)_k=\alpha(\alpha+1)\cdots(\alpha+k-1)$, \quad $k>0$.
\end{quote}
So we have $n\to \infty$,
\begin{quote}
$\frac{\varphi^k\{a_{n-1}\} \varphi^k\{a_{n+1}\} }{ (\varphi^k\{a_{n}\})^2 } < 1$, \quad if $k$ is odd, \\
$\frac{\varphi^k\{a_{n-1}\} \varphi^k\{a_{n+1}\} }{ (\varphi^k\{a_{n}\})^2 } > 1$, \quad if $k$ is even,
\end{quote}
ending the proof of the Theorem.

\begin{ex}
By Theorem ~\ref{Th-2.3}, when $u_n$ are of the following form, the corresponding $\{a_n\}_{n \ge 0}$ satisfies  the  asymptotic $\ell$-log-monotonic inequality.
\[
1 + \frac{1}{n}, \quad  1 + \frac{1}{n \log n}, \quad 1 + \frac{2}{n^2}, \quad 1 + \frac{\log n}{n^2}.
\]
Noting that for $a_n= n!$, we have
\[
u_n = \frac{n+1}{n} = 1 + \frac{1}{n}>0 ,\quad \forall\, n \ge 0.
\]
 Hence $\{ a_n \}_{n \ge 0}$ satisfies the $\ell$-log-monotonic inequality    for $n$ sufficiently large.
\end{ex}

\section{Asymptotically Laguerre inequality of order two}\label{sec3}

In this section, we will give a sufficient condition such that the sequence is asymptotic Laguerre inequality of order two.

\begin{theo} \label{th-3.1}
Let $\{a_n\}_{n \ge 0}$ be a sequence such that \eqref{un} holds.
Assume that $\alpha_m- \alpha_1 \ge 1$. If
\begin{quote}
 $r_1(x)>0$ for $x$ sufficiently large,
\end{quote}
or
\begin{quote}
$\alpha_1<2$, $r_1(x)<0$ for $x$ sufficiently large,
\end{quote}
or
\begin{quote}
$\alpha_1=2$, $r_1(x)<-1$ for $x$ sufficiently large,
\end{quote}
then $\{a_n\}_{n \ge 0}$ satisfies Laguerre inequality of order two asymptotically.
\end{theo}
\proof We only  need to prove if $n\to \infty$, the following  inequality holds
\begin{equation}\label{Lagu}
3a_{n+2}^2-4a_{n+1}a_{n+3}+a_na_{n+4}>0.
\end{equation}
Denote
\[
u_n=\frac{a_{n-1}a_{n+1}}{a_n^2},
\]
dividing the equation ~(\ref{Lagu}) by $a_{n+2}^2$ and by some simple calculation, we need to prove if $n\to \infty$, the following inequalities holds
\begin{equation}
u_{n-1}u_{n+1}^2u_{n+1}-4u_n+3>0.
\end{equation}
Denote
\[
\alpha_1 = \alpha,  \quad
u_n = 1 + \frac{\xi(n)}{n^{\alpha}},
\]
Then we have
\[
\xi(n)= \sum_{i=1}^m  \frac{r_i(\log n)}{n^{\alpha_i-\alpha}} + o\left( \frac{1}{n^{\beta-\alpha}} \right).
\]
Let
\[
f(n) = u_{n-1} u_{n+1}^2 u_{n+1} - 4u_n+3,
\]
then we have
\[
f(n) = \frac{ t_1(n) + t_2(n) + t_3(n) + t_4(n) }{ (n-1)^\alpha  n^{2\alpha} (n+1)^\alpha},
\]
where
\[
t_1(n)=\big[ \xi(n+1) (n-1)^\alpha n^\alpha + \xi(n-1) (n+1)^\alpha  n^\alpha - 2\xi(n) (n-1)^\alpha (n+1)^\alpha \big] n^\alpha,
\]
and
\begin{align*}
t_2(n) & = \xi(n-1) \xi(n+1) n^{2\alpha} + \xi(n)^2 (n-1)^\alpha (n+1)^\alpha\\
       & + 2\xi(n-1) \xi(n) (n+1)^\alpha  n^\alpha + 2\xi(n+1) \xi(n) (n-1)^\alpha  n^\alpha,
\end{align*}
and
\[
t_3(n) = \big[ 2\xi(n-1) \xi(n+1) n^\alpha + \xi(n+1) \xi(n) (n-1)^\alpha + \xi(n-1) \xi(n) (n+1)^\alpha \big] \xi(n),
\]
and
\[
t_4(n) = \xi(n-1) \xi^2(n) \xi(n+1).
\]
By Lemma \ref{lem-2.1} we have
\[
t_1(n) = r_1(\log n)  n^{3\alpha-2} \left( \alpha^2+\alpha+o(1) \right),
\]
and
\[
t_2(n) = 6r_1^2(\log n) n^{2\alpha} \left( 1+o(1) \right),
\]
and
\[
t_3(n) = 4r_1^3(\log n) n^\alpha \left( 1+o(1) \right),
\]
and
\[
t_4(n) = r_1^3(\log n) \left( 1+o(1) \right).
\]
Denote
\[
t(n) = t_1(n) + t_2(n) + t_3(n) + t_4(n).
\]
If $r_1(\log n)>0$, then $n\to \infty$,
\begin{quote}
$t(n)>0$, then $f(n)>0$.
\end{quote}
If $r_1(\log n)<0$, then $n\to \infty$:
\begin{quote}
If $\alpha<2$, the main term of $t(n)$ is $t_2(n)$, we have
\[
f(n) = \frac{6r_1^2(\log n)}{n^{2\alpha}} \left( 1+o(1) \right) > 0.
\]
If $\alpha=2$, the main term of $t(n)$ is $t_1(n)+t_2(n)$, we have
\[
f(n) = \frac{6r_1(\log n) (1+r_1(\log n))}{n^4} \left( 1+o(1) \right).
\]
So if $r_1(\log n)<-1$, we have
\[
f(n)>0.
\]
\end{quote}

\begin{ex}
By Theorem~\ref{th-3.1}, when $u_n$ are of the following form, the corresponding $\{a_n\}_{n \ge 0}$ satisfies  Laguerre inequality of order two asymptotically.
\[
1 \pm\frac{1}{n}, \quad  1 - \frac{1}{n \log n}, \quad 1 - \frac{2}{n^2}, \quad 1 - \frac{\log n}{n^2}.
\]
Noting that for $a_n= n!$, we have
\[
u_n = \frac{n+1}{n} = 1 + \frac{1}{n}.
\]
 Hence $\{ n! \}_{n \ge 0}$ satisfies Laguerre inequality of order two   for $n$ sufficiently large. In fact,
\[
u_{n-1} u_{n+1}^2 u_{n+1} - 4u_n + 3 = \frac{6}{n(n-1)} > 0, \quad \forall\, n \ge 2.
\]
Similarly, if $a_n = \frac{1}{n!}$, we have
\[
u_n  = \frac{n}{n+1} = 1 - \frac{1}{n} + \frac{1}{n^2} + o(n^{-2-\delta}), \quad 0<\delta<1,
\]
 Hence $\{ \frac{1}{n!} \}_{n \ge 0}$ satisfies Laguerre inequality of order two   for $n$ sufficiently large. In fact,
 \[
u_{n-1} u_{n+1}^2 u_{n+1} - 4u_n + 3 = \frac{6}{(n+1)(n+2)} > 0, \quad \forall\, n \ge 2.
\]
\end{ex}

\section{Log-momotonic equality of order three}\label{sec4}
In this section, we will give a method to find the explicit $N$ such that log-momotonic equality of order three holds for
$\{a_n\}_{n\geq N}$ when $\{a_n\}_{n\geq0}$ is P-recursive.

\begin{theo}\label{th-LOG}
Let
\[
u_n = \frac{a_{n-1}a_{n+1}}{a_n^2}.
\]
If there exist an integer $N_1$, an upper bound $f_n$ and a lower bound $g_n$ of $u_n$ such that for all $n \ge N_1$,
\[
0<g_n<u_n<f_n,
\]
and
\[
g_n>1, \quad n\ge N_2,
\]
\[
g_n-f_{n+1}>0, \quad n\ge N_3,
\]
\[
g_{n-1}g_{n+1}-f^2_n>0, \quad n\ge N_4.
\]
Let $N=\max\{N_1,N_2,N_3,N_4\}$, Then  $\{a_n\}_{n\geq N}$ satisfies log-momotonic equality of order three.
\end{theo}
\proof We only  need to prove if $n\ge N$, the following three inequalities holds
\[
\frac{a_{n-1}a_{n+1}}{a_n^2} > 1,
\]
and
\[
\frac{a_n}{a_{n-1}}  \frac{a_{n+2}}{a_{n+1}} < \frac{a_{n+1}^2}{a_n^2},
\]
and
\[
\frac{a_{n-2}a_{n}}{a_{n-1}^2} \frac{a_na_{n+2}}{a_{n+1}^2} > \left( \frac{a_{n-1}a_{n+1}}{a_n^2} \right)^2.
\]
Denote
\[
u_n=\frac{a_{n-1}a_{n+1}}{a_n^2},
\]
By some  simple calculations, we only need to prove $n\ge N$, the following three inequalities holds.
\[
u_n>1,
\]
and
\[
u_n-u_{n+1}>0,
\]
and
\[
u_{n-1}u_{n+1}-u_n^2>0.
\]

If  $n\ge N$, we have
\[
u_n > g_n > 1,
\]
and
\[
u_n - u_{n+1} > g_n - f_{n+1} > 0,
\]
and
\[
u_{n-1} u_{n+1}- u^2_n > g_{n-1} g_{n+1}- f^2_n > 0.
\]

From Theorem ~\ref{th-LOG}, we see the main thing we have to do   is to find a appropriate upper bounds and lower bounds of $u_n$.
We will use the asymptotic estimation of P-recursive sequences given by Birkhoff and Trjitzinsky \cite{Birk} and developed by Wimp and Zeilberger \cite{Wim}. They showed that a P-recursive sequence is asymptotically equal to a linear combination of terms of the form
\begin{equation} \label{Qs}
e^{Q(\rho,n)} s(\rho,n),
\end{equation}
where
\[
Q(\rho,n)=\mu_0 n \log n + \sum_{j=1}^{\rho}\mu_jn^{j/\rho},
\]
and
\[
s(\rho,n)=n^r\sum_{j=0}^{t-1}(\log n)^j\sum_{s=0}^{M-1}b_{sj}n^{-s/\rho},
\]
with $\rho$, $t$, $M$ being positive integers and $\mu_j$, $r$, $b_{sj}$ being complex numbers.

 When $\{a_n\}$ has the asymptotic expression ~(\ref{Qs}) such that $t=1$, we are able to compute the bounds
 $f_n$ and $g_n$ by the  HT algorithm. We implemented the HT algorithm in the package {\tt P}-rec.m
 which is available at \cite{HT}. For more information about the HT algorithm, we refer the author
  to\cite{Houli}.

 \begin{ex}
The central trinomial coefficients $\{b_n\}$, i.e., the coefficients of $x^n$ in expansion of $(1+x+x^2)^n$ satisfies
\[
b_n = \sum_{k=0}^{\lfloor n/2 \rfloor} \binom{n}{2k} \binom{2k}{k}.
\]
Then the sequence $\{b_n\}_{n\geq8}$ satisfies log-momotonic equality of order three.
\end{ex}
\proof
Use Zeilberger's algorithm\cite{Zeilb} we see that $\{b_n\}_{n\geq0}$ satisfies the following recurrence relation
\[
(n+2) b_{n+2} - (2n+3) b_{n+1} - 3(n+1)b_n =  0, \quad n\geq0,
\]
with the initial values
\[
b_0=1,\quad b_1=1.
\]
By the algorithm HT, we have the following estimates:
\[
1 + \frac{1}{2n^2} - \frac{3}{8n^3} + \frac{9}{32n^4} - \frac{355}{256n^5} < \frac{b_{n+1}b_{n-1}}{b_n^2} < 1 + \frac{1}{2n^2} - \frac{3}{8n^3} + \frac{9}{32n^4} + \frac{157}{256n^5}, \quad n\geq12.
\]
Let
\[
u_n = \frac{b_{n+1}b_{n-1}}{b_n^2},
\]
and
\[
g_n = 1 + \frac{1}{2n^2} - \frac{3}{8n^3} + \frac{9}{32n^4} - \frac{355}{265n^5},
\]
and
\[
f_n = 1 + \frac{1}{2n^2} - \frac{3}{8n^3} + \frac{9}{32n^4} + \frac{157}{265n^5}.
\]
Then we have
\[
g_n > 1, \quad n\geq 2.
\]
\[
g_n - f_{n+1} > 0, \quad n\geq 2.
\]
\[
g_{n-1} g_{n+1} -f_n^2 > 0, \quad n\geq 4.
\]
Checking the initial values,  we finally get  $\{b_n\}_{n\geq8}$ satisfies  log-momotonic equality of order three. \qed

 \section{Laguerre inequality of order two}\label{sec5}
 In this section, we will give a method to find the explicit $N$ such that Laguerre inequality of order two holds for
$\{a_n\}_{n\geq N}$ when $\{a_n\}_{n\geq0}$ is P-recursive.

\begin{theo}\label{th-5.1}
Let
\[
u_n = \frac{a_{n-1}a_{n+1}}{a_n^2}.
\]
If there exist an integer $N_1$, an upper bound $f_n$ and a lower bound $g_n$ of $u_n$ such that for all $n \ge N_1$,
\[
0<g_n<u_n<f_n,
\]
and
\[
g_{n-1}  g_n^2  g_{n+1} - 4f_n + 3 > 0, \quad N\geq N_2.
\]
Let $N=\max\{N_1,N_2\}$, Then  $\{a_n\}_{n\geq N}$ satisfies  Laguerre inequality of order two.
\end{theo}
 \proof If  $n\ge N$, we have
\[
u_{n-1} u_{n+1}^2 u_{n+1} - 4u_n+3 > g_{n-1} g_n^2 g_{n+1} - 4 f_n + 3> 0.
\]

\begin{ex}
 Let $M_n$ be the Motzkin numbers given by
\[
(n+4) M_{n+2} - (2n+5) M_{n+1} - 3(n+1) M_n=0, \quad n\geq0,
\]
with the initial values
\[
M_0=1, \quad M_1=1.
\]
Then the sequence $\{\frac{M_n}{n!}\}_{n\geq0}$ satisfies Laguerre inequality of order two.
\end{ex}
\proof
By  the algorithm HT, we derive
\[
1 + \frac{3}{2n^2} - \frac{47}{8n^3} < \frac{M_{n+1}M_{n-1}}{M_n^2} < 1 + \frac{3}{2n^2} - \frac{31}{8n^3}, \quad n\geq228.
\]
Let
\[
a_n = \frac{M_n}{n!}, \quad  u_n = \frac{a_{n+1}a_{n-1}}{a_n^2}.
\]
Then we have
\[
\frac{n}{n+1} \left( 1+\frac{3}{2n^2}-\frac{47}{8n^3} \right) <u_n<\frac{n}{n+1} \left( 1+\frac{3}{2n^2}-\frac{31}{8n^3} \right) , \quad n\geq228.
\]
Let
\[
g_n = \frac{n}{n+1} \left( 1+\frac{3}{2n^2}-\frac{47}{8n^3} \right),
\]
\[
f_n = \frac{n}{n+1} \left( 1+\frac{3}{2n^2}-\frac{31}{8n^3} \right).
\]
Then we have
\[
g_{n-1} g_n^2 g_{n+1}- 4f_n +3 >0,\quad  n\geq 2.
\]

Checking the initial values,  we finally get that $\{\frac{M_n}{n!}\}_{n \ge 0}$ satisfies  Laguerre inequality of order two.

\vskip 0.2cm
\noindent{\bf Acknowledgments.}
 We would like to thank the referees for their helpful comments.

\end{document}